\documentclass[twoside,leqno]{article}

\usepackage[letterpaper]{geometry}

\usepackage{ltexpprt}

\newcommand{\qobisubfigure}[2]{\begin{tabular}[t]{@{}c@{}}#2\\#1\end{tabular}}
\usepackage{booktabs}
\usepackage{xcolor}

\usepackage{amsfonts}
\usepackage{amsmath}
\usepackage{fdsymbol}
\usepackage{graphicx}
\usepackage{accents}
\usepackage{textcomp}

\usepackage{hyperref}
\definecolor{RED}{rgb}{1,0,0}
\definecolor{darkred}{rgb}{0.85,0,0}
\definecolor{darkgreen}{rgb}{0,0.6,0}
\definecolor{darkblue}{rgb}{0,0,0.85}
\hypersetup{hidelinks, breaklinks, colorlinks,
  linkcolor={darkblue}, citecolor={darkgreen}, urlcolor={darkblue}}

\graphicspath{{figures/}}

\newcommand{\citet}[1]{\cite{#1}}
\newcommand{\citep}[1]{\cite{#1}}

\DeclareUnicodeCharacter{03BB}{\lambda} 
\DeclareUnicodeCharacter{FE52}{\;.\;} 
\DeclareUnicodeCharacter{03C0}{\pi} 
\DeclareUnicodeCharacter{211D}{\Re} 
\DeclareUnicodeCharacter{2192}{\rightarrow} 
\DeclareUnicodeCharacter{21A6}{\mapsto} 
\DeclareUnicodeCharacter{21A4}{\leftmapsto} 
\DeclareUnicodeCharacter{03B1}{\alpha} 
\DeclareUnicodeCharacter{1D7CE}{\mm{0}} 
\DeclareUnicodeCharacter{1D6AB}{\Delta} 
\DeclareUnicodeCharacter{1D400}{\mm{A}} 
\DeclareUnicodeCharacter{1D401}{\mm{B}} 
\DeclareUnicodeCharacter{1D408}{\mm{I}} 
\DeclareUnicodeCharacter{1D409}{\mm{J}} 
\DeclareUnicodeCharacter{1D4A5}{\mathcal{J}} 
\DeclareUnicodeCharacter{2148}{\mathrm{i}} 
\DeclareUnicodeCharacter{1D53B}{\DD} 
\DeclareUnicodeCharacter{1D49F}{\D} 
\DeclareUnicodeCharacter{27F9}{\reducesto} 
\DeclareUnicodeCharacter{2080}{_0} 
\DeclareUnicodeCharacter{2081}{_1} 
\DeclareUnicodeCharacter{2082}{_2} 
\DeclareUnicodeCharacter{2070}{^0} 
\DeclareUnicodeCharacter{2071}{^1} 
\DeclareUnicodeCharacter{00B2}{^2} 
\DeclareUnicodeCharacter{00B3}{^3} 
\DeclareUnicodeCharacter{FF5C}{\mid} 
\DeclareUnicodeCharacter{03B5}{\ve} 
\DeclareUnicodeCharacter{25CB}{\circ} 
\DeclareUnicodeCharacter{2219}{\cdot} 
\DeclareUnicodeCharacter{2026}{\ldots} 
\DeclareUnicodeCharacter{22EF}{\cdots} 
\DeclareUnicodeCharacter{00D7}{\times} 
\DeclareUnicodeCharacter{225D}{\define} 

\renewcommand{\Re}{\mathbb{R}}
\newcommand{\vv}[1]{\mathbf{#1}}	
\newcommand{\mm}[1]{\mathbf{#1}}	

\newcommand{\fv}[1]{\dot{#1}}
\newcommand{\rv}[1]{\bar{#1}}
\newcommand{\fiv}[1]{\dot{#1}^*}
\newcommand{\riv}[1]{\bar{#1}^*}
\newcommand{\inv}[1]{#1^{-1}}
\newcommand{\tran}[1]{#1^{\mathsf{T}}}
\newcommand{\itran}[1]{#1^{-\mathsf{T}}}
\newcommand{\myfrac}[2]{\frac{\raisebox{-0.5ex}{\rule{0em}{0ex}}#1}{\text{\raisebox{0.5ex}{$#2$}}}}
\newcommand{\define}{\stackrel{\triangle}{=}}

\newcommand{\ddt}[0]{\frac{\mathrm{d}}{\mathrm{d}t}}

\newcommand{\J}{{\overrightarrow{𝒥}}}
\newcommand{\Jforward}{\J}
\newcommand{\Jreverse}{\text{\reflectbox{$\J$}}}
\newcommand{\Jforwardinverse}{\text{\raisebox{-5pt}{\raisebox{\depth}{\rotatebox{180}{\reflectbox{$\J$}}}}}}
\newcommand{\Jreverseinverse}{\text{\raisebox{-5pt}{\raisebox{\depth}{\rotatebox{180}{$\J$}}}}}

\usepackage{xspace}
\def\onedot{\ifx\@let@token.\else.\null\fi\xspace}
\newcommand{\eg}{\emph{e.g.},}

\newcommand{\ie}{\emph{i.e.},}

\newcommand{\etc}{{\emph{etc}}\onedot}

\newcommand{\vs}{\emph{vs.}}

\begin{document}

\newcommand\relatedversion{}
\renewcommand\relatedversion{\thanks{The full version of the paper can be accessed at \protect\url{https://arxiv.org/abs/1902.09310}}} 

\title{\Large Automatic Differentiation: Inverse Accumulation Mode}
\author{\href{http://barak.pearlmutter.net}{Barak A. Pearlmutter}\thanks{Department of Computer Science, Maynooth University,
    Maynooth, Co.\ Kildare, Ireland}
  \and
  \href{http://www.ece.purdue.edu/~qobi/}{Jeffrey Mark Siskind}\thanks{Elmore Family School of Electrical and Computer Engineering,
    Purdue University, West Lafayette, IN 47907-2035, USA}
}



\date{}

\maketitle


\fancyfoot[R]{\scriptsize{Copyright \textcopyright\ 2024\\
Copyright for this paper is retained by authors.}}

\begin{abstract} \small\baselineskip=9pt
  We show that, under certain circumstances, it is possible to
  automatically compute Jacobian-inverse-vector and
  Jacobian-inverse-transpose-vector products about as efficiently as
  Jacobian-vector and Jacobian-transpose-vector products.
  The key insight is to notice that the Jacobian corresponding to the use of
  one basis function is of a form whose sparsity is invariant to inversion.
  The main restriction of the method is a constraint on the number of active
  variables, which suggests a variety of techniques or generalization to allow
  the constraint to be enforced or relaxed.
  This technique has the potential to allow the efficient direct
  calculation of Newton steps as well as other numeric calculations of interest.
\end{abstract}

\section{Inverse AD: The Dream}

Automatic Differentiation (AD) is the mechanical transformation of computer
programs to calculate derivatives of interest, with useful complexity
guarantees.
The two most important ``modes'' of AD are forward and reverse, which
access the Jacobian (the matrix of derivatives of each output of the
computation with respect to each input) by multiplication, or
transpose-multiplication, with a vector.
Here we consider first-order numeric computations, where inputs and
outputs are vectors of reals.
Given the primal computation
\begin{math}
  y = f(x)
\end{math}
with $f : ℝ^m → ℝ^n$ and therefore $x\in ℝ^m$ and
$y\in ℝ^n$, we use $𝐉_{f(x)} \in ℝ^{n×m}$ for the Jacobian of
$f$ at $x$, whose ${(i,j)}^{\text{th}}$ element is $\partial
f_i(x)/\partial x_j$.
Forward and Reverse AD compute
\begin{align} \label{eq:AD}
  \fv{y} &= 𝐉_{f(x)} \: \fv{x}
  &
  \text{and}
  &&
  \rv{x} &= \tran{𝐉_{f(x)}} \: \rv{y}
\end{align}
respectively, with $\fv{\cdot}$ and $\rv{\cdot}$ denoting tangents and
cotangents.
Our objective here is to find an efficient way to solve for the
starred vectors on the right-hand sides in each of
\begin{subequations} \label{eq:solve}
\begin{align}
  \riv{x} &= \tran{𝐉_{f(x)}} \: \riv{y}
  \label{eq:forwardInverse}
  \\
  \fiv{y} &= 𝐉_{f(x)} \: \fiv{x}
  \label{eq:reverseInverse}
\end{align}
\end{subequations}
If this can be done efficiently, it would allow efficient Newton steps (where
$f$ is a gradient calculation, say) and other sorts of second-order
optimization.
For this to be well posed it is necessary for $𝐉_{f(x)}$ to be
invertible, so $n=m$, and
as we shall see, further restrictions on the form of $f$ will be required.
Inverse Jacobians can be used to find roots of systems of equations.
Inverse Hessians (which can be computed with Inverse AD over
traditional AD) can be used for second order optimization, and for this reason are the topic of intensive research in the optimization community \citep{bullins-etal-2021a}.

\section{Inverse AD: The Reality}

Let us review Forward and Reverse AD.
Since we are evaluating $f$ at a point $x$, we consider control flow resolved
and represent the computation as a data flow graph: a DAG whose edges hold
reals and whose vertices represent numeric basis functions: things like $+$,
$-$, $\times$, $\div$, $\sqrt{\cdot}$, $\log$, $\exp$, $\sin$, $\cos$,
$\tan^{-1}$, \etc, that are primitives, intrinsics, or library functions in
a typical programming language.
There are $n$ edges entering from the inputs $x₁, …, x_n$, and $n$
exiting to $y₁, …, y_n$.
If we topologically sort the data flow graph, and cut it before and
after each vertex, we see that the computation
proceeds through a sequence of $T+1$ machine states, $\vv{x}₀, …,
\vv{x}_T$, where the initial and final states are the input and output of the
computation, $\vv{x}₀ = x$ and $y = \vv{x}_T$.
We can denote the transition function from one machine state to the next by
$\vv{x}_t = f_t(\vv{x}_{t-1})$ and the Jacobian of $f_t$ at $\vv{x}_{t-1}$
by $𝐉_t$, keeping in mind that $f_t$ involves applying a single numeric
basis function to some elements of $\vv{x}_{t-1}$ and putting the result in
some elements of $\vv{x}_t$, copying the other elements unchanged.

Since $f = f_T ○ f_{T-1} ○ ⋯ ○ f₂ ○ f₁$, the Jacobian
matrix is a product, $𝐉_{f(x)} = 𝐉_T \: 𝐉_{T-1} ⋯ 𝐉₂ \:
𝐉₁$, and Forward and Reverse AD amount to appropriate associativity:
\begin{subequations} \label{eq:farAD}
\begin{align}
  \fv{y} &= 𝐉_{f(x)} \, \fv{x} = 𝐉_T (𝐉_{T-1} ⋯ (𝐉₂ (𝐉₁ \: \fv{x}))⋯)
  \\
  \rv{x} &= \tran{𝐉}_{f(x)} \, \rv{y} = \tran{𝐉₁} (\tran{𝐉₂} ⋯ (\tran{𝐉_{T-1}} (\tran{𝐉_T} \: \rv{y})) ⋯)
\end{align}
\end{subequations}

\begin{figure*}[t!]
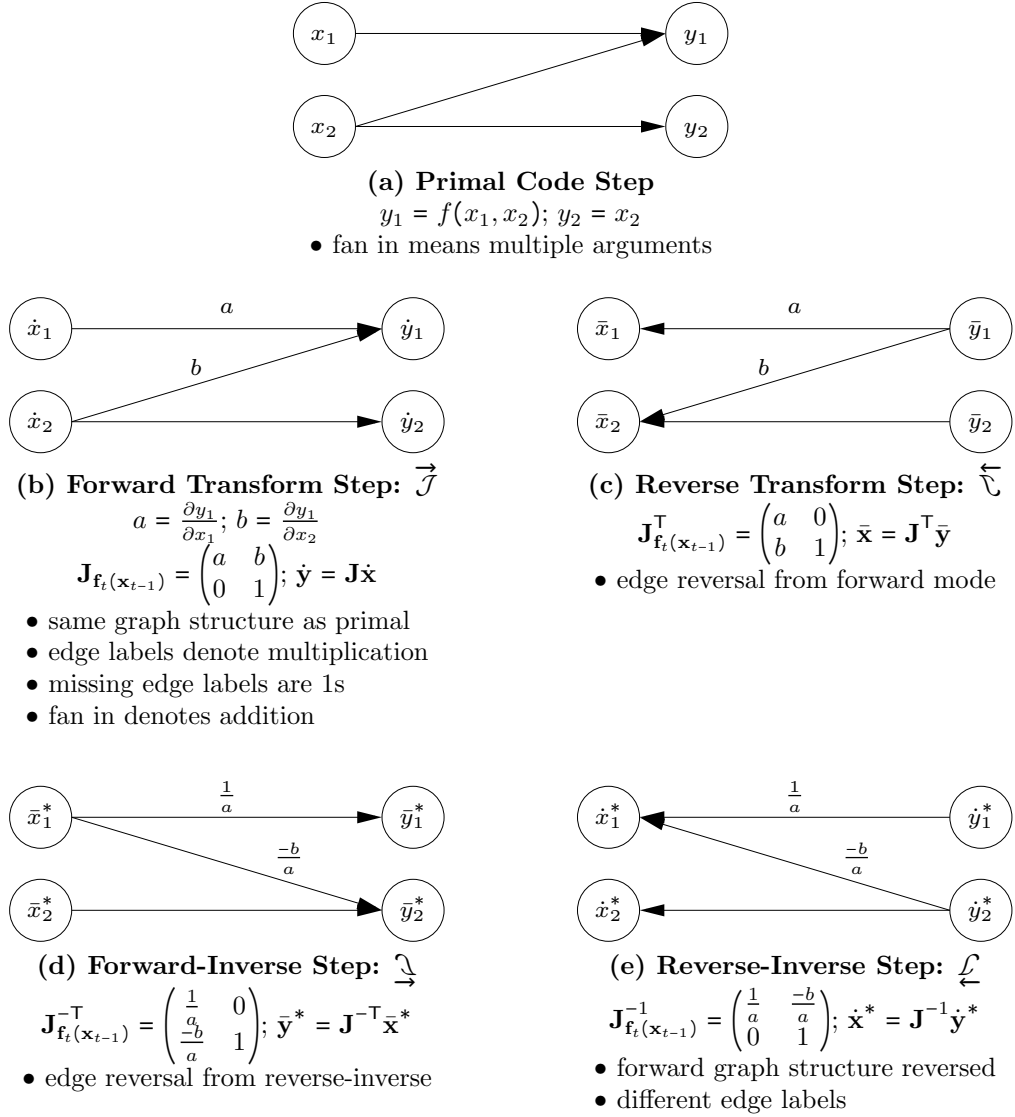

  \centering
  \begin{tabular}{c@{\hspace{5em}}c}
    \multicolumn{2}{c}{
    \qobisubfigure{\textbf{(a) Primal Code Step}\\
      $y_1 = f(x_1,x_2)$; $y_2 = x_2$\\
      $\bullet$ fan in means multiple arguments
    }{\scalebox{1.25}{\input{incl/primal-essence-binary}}}
    }
    \\[10ex]
    \qobisubfigure{\textbf{(b) Forward Transform Step: $\Jforward$}\\
      $a=\frac{\partial y_1}{\partial x_1}$;
      $b=\frac{\partial y_1}{\partial x_2}$
      \\
      $\mathbf{J}_{\mathbf{f}_t(\mathbf{x}_{t-1})}=
      \begin{pmatrix}
        a & b\\
        0 & 1\\
      \end{pmatrix}$;
      $\fv{\mathbf{y}}=\mathbf{J}\fv{\mathbf{x}}$
      \\
      \begin{tabular}{l}
      $\bullet$ same graph structure as primal\\
      $\bullet$ edge labels denote multiplication\\
      $\bullet$ missing edge labels are 1s\\
      $\bullet$ fan in denotes addition
      \end{tabular}
    }{\scalebox{1.25}{\input{incl/forward-essence-binary}}}
    &
    \qobisubfigure{\textbf{(c) Reverse Transform Step: $\Jreverse$}\\
      $\tran{\mathbf{J}}_{\mathbf{f}_t(\mathbf{x}_{t-1})}=
      \begin{pmatrix}
        a & 0\\
        b & 1\\
      \end{pmatrix}$;
      $\rv{\mathbf{x}}=\tran{\mathbf{J}}\rv{\mathbf{y}}$
      \\
      $\bullet$ edge reversal from forward mode
    }{\scalebox{1.25}{\input{incl/reverse-essence-binary}}}
    \\[27ex]
    \qobisubfigure{\textbf{(d) Forward-Inverse Step: $\Jforwardinverse$}\\
      $\itran{\mathbf{J}}_{\mathbf{f}_t(\mathbf{x}_{t-1})} =
      \begin{pmatrix}
        \frac{1}{a} & 0\\
        \frac{-b}{a} & 1\\
      \end{pmatrix}$;
      $\riv{\mathbf{y}}=\itran{\mathbf{J}}\riv{\mathbf{x}}$
      \\
      $\bullet$ edge reversal from reverse-inverse
    }{\scalebox{1.25}{\input{incl/forward-inverse-essence-binary}}}
    &
    \qobisubfigure{\textbf{(e) Reverse-Inverse Step: $\Jreverseinverse$}\\
      $\inv{\mathbf{J}}_{\mathbf{f}_t(\mathbf{x}_{t-1})} =
      \begin{pmatrix}
        \frac{1}{a} & \frac{-b}{a}\\
        0 & 1\\
      \end{pmatrix}$;
      $\fiv{\mathbf{x}}=\inv{\mathbf{J}}\fiv{\mathbf{y}}$
      \\
      \begin{tabular}{l}
      $\bullet$ forward graph structure reversed\\
      $\bullet$ different edge labels
      \end{tabular}
    }{\scalebox{1.25}{\input{incl/reverse-inverse-essence-binary}}}
  \end{tabular}
  \caption{Graphical representation of transformation of computation graph of
    binary atomic program step, for all four AD modes discussed.
    These are formulated for scalar inputs and outputs.
    In the case where the first input/output is a vector of length $l$ and the
    second input is a vector of length $k-l$, one simply replaces $a$ with $𝐀$,
    $b$ with $𝐁$, $\frac{1}{a}$ with $\inv{𝐀}$, $\frac{-b}{a}$ with
    $- \inv{𝐀} 𝐁$, $0$ with $𝟎$, and $1$ with $𝐈$.}
  \label{fig:graphAtomicXforms}
\end{figure*}

\begin{figure*}
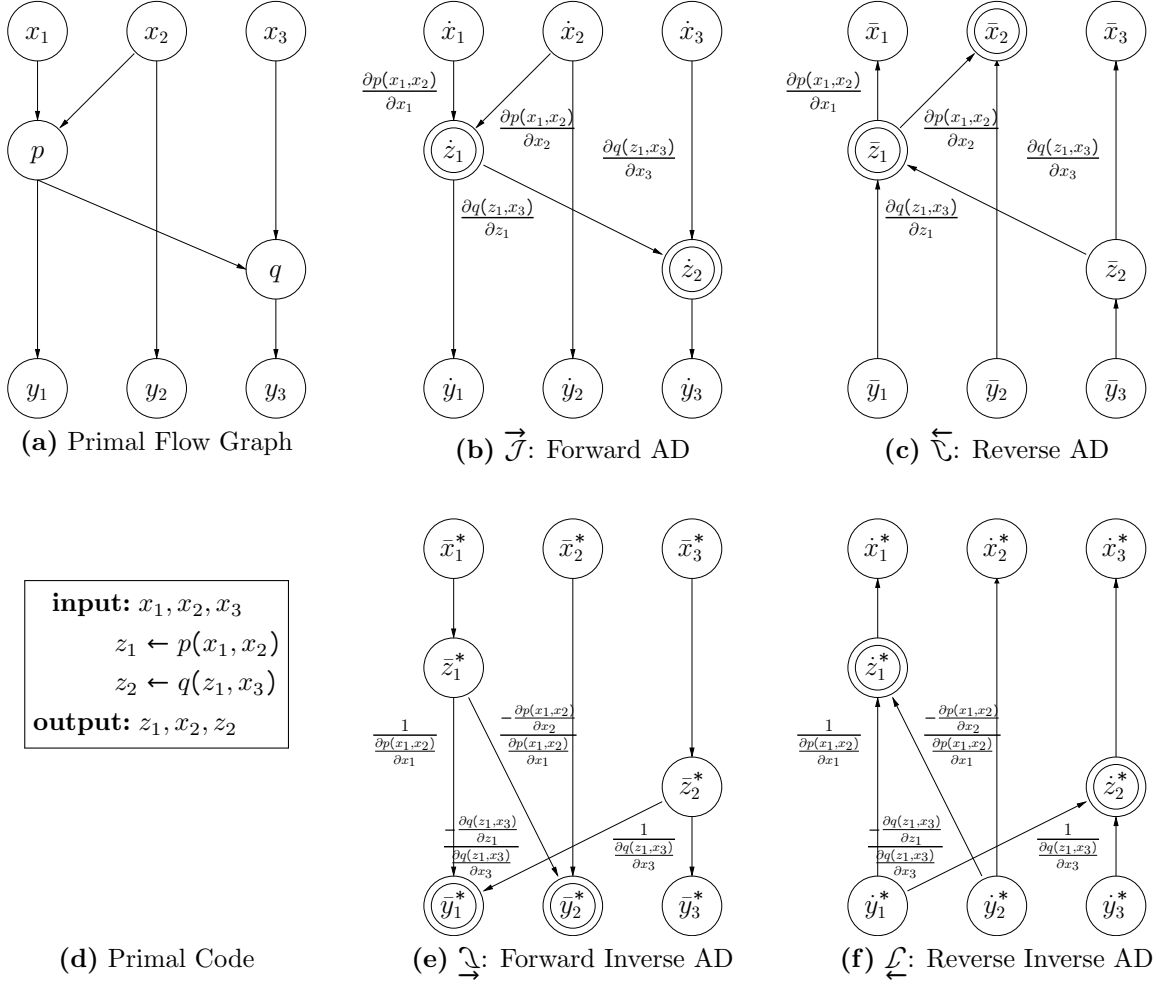

  \centering
  \newcommand{\diagramscale}{0.6}
  \begin{tabular}{c@{\hspace{4em}}c@{\hspace{4em}}c}
      \qobisubfigure{\textbf{(a)} Primal Flow Graph}{\scalebox{\diagramscale}{\input{incl/dag-primal}}}&
      \qobisubfigure{\textbf{(b)} $\Jforward$: Forward AD}{\scalebox{\diagramscale}{\input{incl/dag-forward}}}&
      \qobisubfigure{\textbf{(c)} $\Jreverse$: Reverse AD}{\scalebox{\diagramscale}{\input{incl/dag-reverse}}}\\[8ex]
      \qobisubfigure{\textbf{(d)} Primal Code}{\raisebox{100pt}{\fbox{\(
          \begin{aligned}
            \textbf{input:}\;&x₁, x₂, x_3\\
            z₁&\leftarrow p(x₁, x₂)\\
            z₂&\leftarrow q(z₁, x_3)\\
            \textbf{output:}\;&z₁, x₂, z₂
          \end{aligned}\)}}}&
      \qobisubfigure{\textbf{(e)} $\Jforwardinverse$: Forward Inverse AD}{\scalebox{\diagramscale}{\input{incl/dag-forward-inverse}}}&
      \qobisubfigure{\textbf{(f)} $\Jreverseinverse$: Reverse Inverse AD}{\scalebox{\diagramscale}{\input{incl/dag-reverse-inverse}}}
  \end{tabular}
  \caption{Illustration of all four AD modes for the straight-line code in~(d).
    This corresponds to the data flow graph~(a).
    The intent is that there are three registers, $r₁$, $r₂$, and~$r_3$,
    illustrated by the three columns in~(a) from left to right.
    These are initialized with~$x₁$, $x₂$, and~$x_3$ respectively.
    Since~$r₁$ is not used after the first line of code, it is overwritten
    with~$z₁$.
    Since~$r_3$ is not used after the second line of code, it is overwritten
    with~$z₂$.
    Forward mode and reverse mode are shown in~(b) and~(c) respectively.
    In these graphs, addition occurs whenever there is fan in to a vertex (the
    circled vertices) and labels on edges denote multiplication by the
    indicated coefficient.
    Reverse mode is derived from forward mode by edge reversal, which can
    change which vertices perform addition due to fan in.
    Forward inverse mode and reverse inverse mode are shown in~(e) and~(f)
    respectively.
    These have the same vertices as forward mode and reverse mode but
    different edges and edge labels, which changes which vertices perform
    addition due to fan in.
    Again, forward inverse mode is derived from reverse inverse mode by edge
    reversal.}
  \label{fig:example}
\end{figure*}

\begin{figure*}[t!]
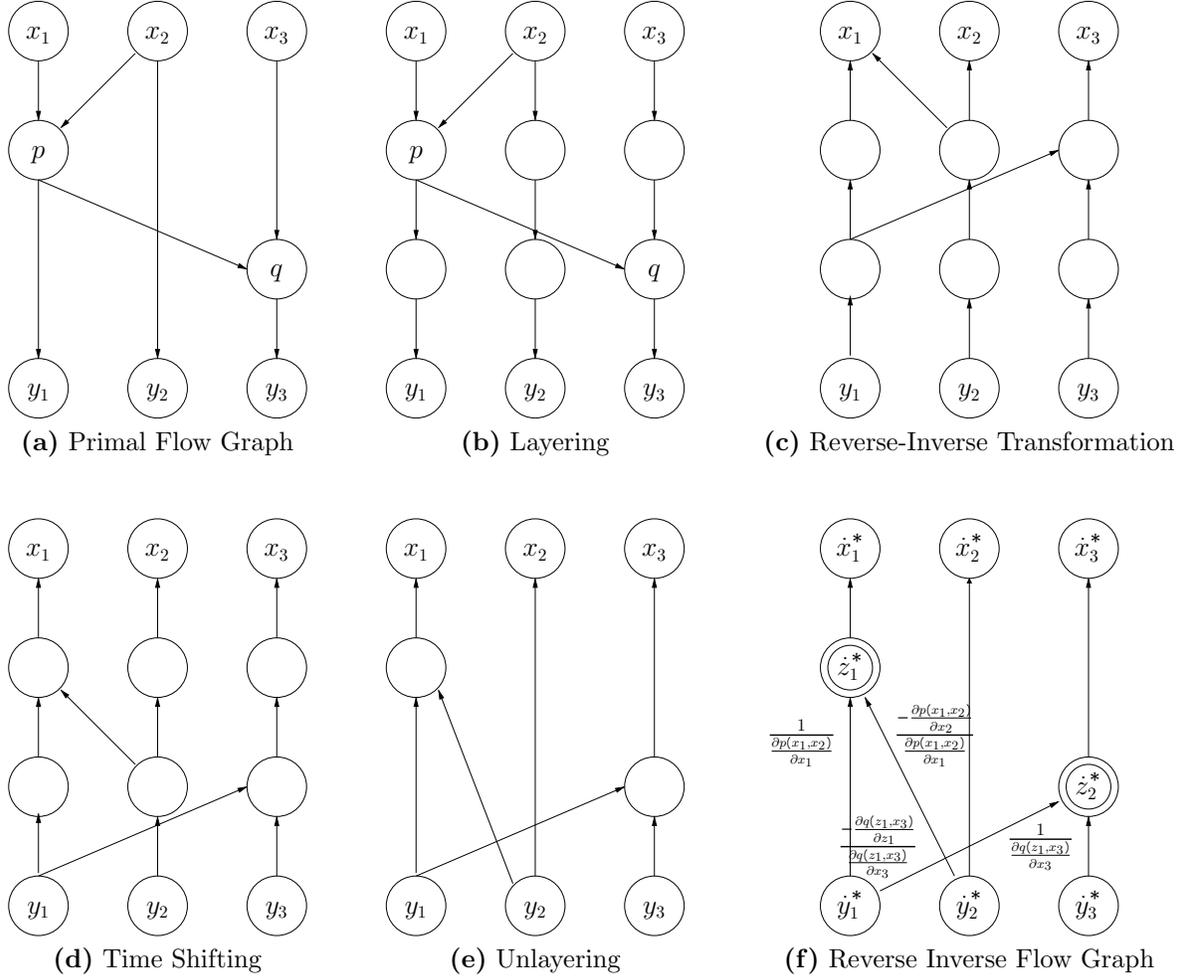

  \centering
  \newcommand{\diagramscale}{0.6}
  \begin{tabular}{c@{\hspace{3em}}c@{\hspace{3em}}c}
      \qobisubfigure{\textbf{(a)} Primal Flow Graph}{\scalebox{\diagramscale}{\input{incl/dag-primal}}}&
      \qobisubfigure{\textbf{(b)} Layering}{\scalebox{\diagramscale}{\input{incl/dag-primal-layered}}}&
      \qobisubfigure{\textbf{(c)} Reverse-Inverse Transformation}{\scalebox{\diagramscale}{\input{incl/dag-reverse-inverse-layered}}}\\[8ex]
      \qobisubfigure{\textbf{(d)} Time Shifting}{\scalebox{\diagramscale}{\input{incl/dag-reverse-inverse-layered2}}}&
      \qobisubfigure{\textbf{(e)} Unlayering}{\scalebox{\diagramscale}{\input{incl/dag-reverse-inverse-unlayered}}}&
      \qobisubfigure{\textbf{(f)} Reverse Inverse Flow Graph}{\scalebox{\diagramscale}{\input{incl/dag-reverse-inverse}}}
  \end{tabular}
  \caption{Illustration of the derivation of Fig.~\ref{fig:example}(f) from
    Fig.~\ref{fig:example}(a).
    Panel~(a) corresponds to Fig.~\ref{fig:example}(a).
    Panel~(b) corresponds to construction of a layered flow graph \citep{naumann-2024a}
    by carrying live variables forward.
    Panel~(c) corresponds to applying the transformation of
    Fig.~\ref{fig:graphAtomicXforms}(e).
    Panel~(d) corresponds to shifting each operation one time step earlier to
    eliminate the noop in the first time step.
    Panel~(e) corresponds to removing the layering.
    Panel~(f) corresponds to Fig.~\ref{fig:example}(f).
    Fig.~\ref{fig:example}(b,c) are derived from Fig.~\ref{fig:example}(a) using
    standard AD methods.
    Fig.~\ref{fig:example}(c) is derived from Fig.~\ref{fig:example}(b) using
    edge reversal.
    Fig.~\ref{fig:example}(e) is derived from Fig.~\ref{fig:example}(f) using
    edge reversal.}
  \label{fig:derived}
\end{figure*}

Solving \eqref{eq:solve} in the form of \eqref{eq:farAD} while
assuming each $𝐉_t$ is invertible is the basic idea of
\emph{Forward Inverse Accumulation} and \emph{Reverse Inverse
  Accumulation:}
\begin{subequations} \label{eq:inverseAD}
\begin{align}
  \riv{y} &= \itran{𝐉_T} (\itran{𝐉_{T-1}}⋯ (\itran{𝐉₂} (\itran{𝐉₁} \: \riv{x})) ⋯)
  \\
  \fiv{x} &= \inv{𝐉₁} (\inv{𝐉₂} ⋯ (\inv{𝐉_{T-1}} (\inv{𝐉_T} \: \fiv{y})) ⋯)
\end{align}
\end{subequations}
where $\itran{\mathbf{M}}=\tran{(\inv{\mathbf{M}})}$.
These will be practical if the matrix-vector products $\inv{𝐉_t}
\fiv{y}$ and $\itran{𝐉_t} \riv{x}$ can be calculated efficiently.
Assuming the computation of $f$ is constant-width, so $\vv{x}_t \in ℝ^n$, and
its local linearization is invertible, then each $f_t$ must write its result to
a slot where one of the inputs to the invoked basis function was stored,
yielding Jacobians of the form
\begin{subequations}\label{eq:primitives}
\begin{align}
    \begin{array}[t]{c}
      \begin{array}{r|ccccccc|}
        \multicolumn{1}{r}{}
        & &      & &R_t       & & &\multicolumn{1}{c}{}\\
        \multicolumn{1}{r}{}
        & &      & &\downarrow& & &\multicolumn{1}{c}{}\\
        \cline{2-8}
        &1&      & &          & & &\\
        & &\ddots& &          & & &\\
        & &      &1&          & & &\\
R_t→&& & &    a     & & &\\
        & &      & &          &1& &\\
        & &      & &          & &\ddots&\\
        & &      & &          & &      &1\\
        \cline{2-8}
      \end{array}\\
      \vspace*{-1ex}\\
      \hspace{3em} a=\frac{\partial g(\vv{x}_{t-1}[R_t])}{\partial \vv{x}_{t-1}[R_t]}
    \end{array}\\
    \intertext{for unary basis functions $g$, that read and
write to variable/slot $R_t$, and}
    \begin{array}[t]{c}
      \begin{array}{r|ccccccc|}
        \multicolumn{1}{r}{}
        & &      & &R_t       & &S_t       &\multicolumn{1}{c}{}\\
        \multicolumn{1}{r}{}
        & &      & &\downarrow& &\downarrow&\multicolumn{1}{c}{}\\
        \cline{2-8}
        &1&      & &          & &          & \\
        & &\ddots& &          & &          & \\
        & &      &1&          & &          & \\
        R_t→
        & &      & &  a       & & b       & \\
        & &      & &          &1&         & \\
        & &      & &          & &\ddots   & \\
        & &      & &          & &         & 1\\
        \cline{2-8}
      \end{array}\\
      \vspace*{-1ex}\\
      \hspace{3.5em} a=\frac{\partial h(\vv{x}_{t-1}[R_t],\vv{x}_{t-1}[S_t])}{\partial \vv{x}_{t-1}[R_t]}\\[2.5ex]
      \hspace{3.5em} b=\frac{\partial h(\vv{x}_{t-1}[R_t],\vv{x}_{t-1}[S_t])}{\partial \vv{x}_{t-1}[S_t]}\\
    \end{array}
\end{align}
\end{subequations}
for binary basis functions $h$ that read from variables/slots $R_t$ and $S_t$
and write to variable/slot $R_t$
We now note that these can be trivially inverted.
If we consider only variables involved in the basis function being invoked,
and reorder them so the output values are first, a basis function with $k$
inputs and a scalar output results in
\begin{subequations} \label{eq:blocks}
\begin{align}  \label{eq:blockrow}
  𝐉_t &= \left( \begin{array}{c|ccc}
    a & b₁ & ⋯ & b_{k-1} \\\hline
    \rule{0em}{2.5ex} 𝟎 & \multicolumn{3}{c}{𝐈}
  \end{array}\right)
  \\[2ex]
  \inv{𝐉_t} &= \left( \begin{array}{c|ccc}
    \myfrac{1}{a} & -\myfrac{b₁}{a} & ⋯ & -\myfrac{b_{k-1}}{a} \\[1pt]\hline
    \rule{0em}{2.5ex} 𝟎 & \multicolumn{3}{c}{𝐈}
  \end{array}\right)
\end{align}
We can generalize from scalar to $l$ outputs, giving the form
\begin{align} \label{eq:blockblock}
  𝐉_t &= \begin{pmatrix} 𝐀 & 𝐁 \\ 𝟎 & 𝐈 \end{pmatrix}
  \\
  \inv{𝐉_t} &= \begin{pmatrix} \inv{𝐀} & - \inv{𝐀} 𝐁 \\ 𝟎 & 𝐈 \end{pmatrix}
\end{align}
where $𝐀:l×l$ and $𝐁:l×(k-l)$.
\end{subequations}
Although $𝐉_t$ is not structurally symmetric, $\inv{𝐉_t}$ has the same
structural sparsity as $𝐉_t$.
And although the amount of arithmetic is the same as for conventional Forward and
Reverse modes, these are transposed, so Forward Inverse Mode writes to the
derivative-related quantities associated with \emph{all} involved variables of
each basis function invocation, while Reverse Inverse Mode writes only to the
quantities associated with slots \emph{written to} in the primal computation
of each basis function.
Fig.~\ref{fig:graphAtomicXforms} shows how all four AD modes transform atomic
portions of a computation graph.
Fig.~\ref{fig:example} illustrates all four AD modes on a simple program.
Fig.~\ref{fig:derived} illustrates how Fig.~\ref{fig:example}(bcef) are derived
from Fig.~\ref{fig:example}(a) with the transformations of
Fig.~\ref{fig:graphAtomicXforms} together with the layering technique of
\citet{naumann-2024a}.

Traditional forward mode can be computed without saving intermediate values
from the primal, as the primal and tangent can be computed in tandem.
However, traditional reverse mode requires a tape to save the intermediate
values from the primal computed during the forward sweep for use in reverse
order during the reverse sweep.
Analogously, forward inverse mode does not require a tape while reverse inverse
mode does.

It is apparent from \eqref{eq:blockblock} that any computation step that takes
$k$ inputs and produces $l$ outputs will have an invertible Jacobian if $𝐀:l×l$
is nonsingular.
In particular, when $l=1$, that will be when $a\not=0$.

\section{Constant Width Graph} \label{sec:constWidth}

The limitation of the methods proposed above is that they require the
computation to be constant width.
What that means is that when the overall function is $ℝ^n→ℝ^n$, there are
precisely $n$ live active variables (active in the AD sense) at each
intermediate point between computation steps.
This requires that every output value of a computation step
overwrite some input value.
Not all programs have this property.
Some programs have varying numbers of live active variables as the computation
proceeds, \ie\ temporary variables.
If ever the number of live active variables is less than $n$, the
Jacobian of the computation is necessarily singular.

\begin{figure*}
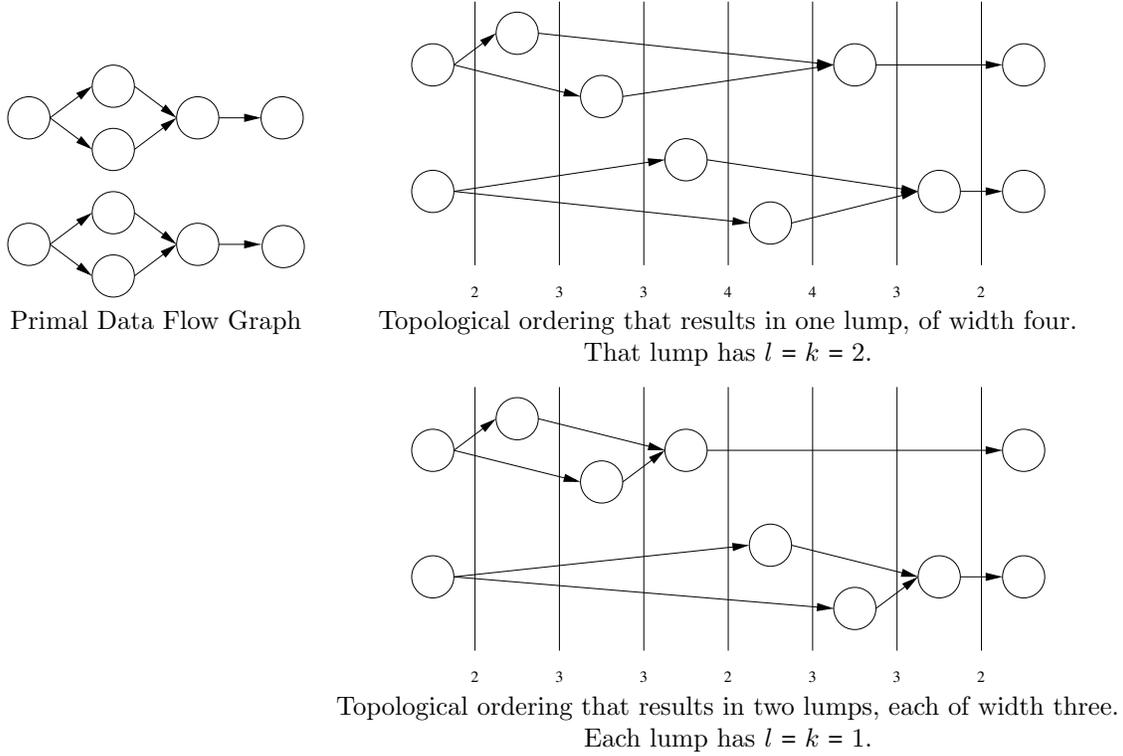

  \centering
  \begin{tabular}{cc}
  \qobisubfigure{Primal Data Flow Graph}{\scalebox{0.85}{\input{incl/graph1}}}
  &
  \qobisubfigure{Topological ordering that
    results in one lump, of width four.\\
    That lump has $l=k=2$.
  }{\scalebox{0.85}{\input{incl/graph2}}}
  \\[7ex]
  &
  \qobisubfigure{Topological ordering that
    results in two lumps, each of width three.\\
    Each lump has $l=k=1$.
  }{\scalebox{0.85}{\input{incl/graph3}}}
  \end{tabular}
  \caption{Illustration of lumpification's dependence on how a total ordering
    is imposed on the partially-ordered data-flow graph.}
  \label{fig:lumps}
\end{figure*}

But if ever the number of live active variables is greater than $n$, the
computation can be partitioned into ``lumps'' where the number of live active
variables between lumps is $n$.
Each lump can be
treated as a macro step and processed according to \eqref{eq:blocks}.
The question then reduces to
how to best partition a computation into lumps, and whether when doing so $l$ is
small.
This process can be viewed as performing a topological sort of the computation
graph, and then breaking that fully-ordered computation up into ``lumps'' at
points where there are exactly $n$ live variables.
But there are many possible topological sorts, which may break the computation
into different numbers of lumps of different sizes, as shown in
Fig.~\ref{fig:lumps}.
We conjecture that a ``greedy'' algorithm will perform optimally here, but
proving this conjecture is left for future work.
We have implemented a system using preliminary answers to this question, which
will be exhibited when its efficiency properties have been further explored.

\section{Notation}

We have a bit of a combinatorial explosion of AD modes on our hands: forward
\vs\ reverse, and noninverse \vs\ inverse, yielding four modes.
To reason about these more easily, and facilitate their inclusion as
first-class operators in programming languages, we propose notation which is
short and mnemonic.
Using a horizontal flip to indicate forward \vs\ reverse and a vertical one to
indicate inverse, with a calligraphic J for Jacobian surmounted by an arrow to
indicate flow, we have:
\begin{center}
    \begin{tabular}{rc@{\qquad}c}\toprule
      & forward & reverse\\
      \cmidrule{2-3}
      noninverted &
      $\Jforward\;f\;x\;\fv{x} ≝ 𝐉_{f(x)} \: \fv{x}$ &
      $\Jreverse\;f\;x\;\rv{y} ≝ \tran{𝐉_{f(x)}} \: \rv{y}$\\
      inverted &
      $\Jforwardinverse\;f\;x\;\riv{x} ≝ \itran{𝐉_{f(x)}} \: \riv{x}$ &
      $\Jreverseinverse\;f\;x\;\fiv{y} ≝ \inv{𝐉_{f(x)}} \: \fiv{y}$\\[3pt]
      \bottomrule
    \end{tabular}
\end{center}
which obey a number of algebraic invariants:
\begin{subequations} \label{eq:invariants}
\begin{gather}
  (\Jreverse\;f\;x\;\rv{y})∙\fv{x} = \rv{y}∙(\Jforward\;f\;x\;\fv{x})
  \\
  \riv{x}∙(\Jreverseinverse\;f\;x\;\fiv{y}) = (\Jforwardinverse\;f\;x\;\riv{x})∙\fiv{y}
  \\
  \Jforward\;f\;x○\Jreverseinverse\;f\;x = 
  \Jforwardinverse\;f\;x ○ \Jreverse\;f\;x = \textbf{id}
  \\
  \Jreverse\;f\;x○\Jforwardinverse\;f\;x = 
  \Jreverseinverse\;f\;x○\Jforward\;f\;x = \textbf{id}
\end{gather}
\end{subequations}

\section{Inverse AD of an ODE}

Consider a primal computation $x₀ ↦ x_T$ where $x(T₀) = x₀$ is the initial
condition of an ODE
\begin{equation} \label{eq:ode}
  \ddt x = g(x)
\end{equation}
and $x_T = x(T₁)$ is its final condition, the result of integrating
\eqref{eq:ode} from $T₀$ to~$T₁$.
(If we wish to give $g$ side parameters we can concatenate them onto $x$ and
extend $g$ to give them zero derivatives, thus incorporating them into the
current treatment without loss of generality.
Alternatively, these can be treated as constants rather than active variables,
so as non-active variables they do not enter into the constant-width
calculation.)

We will discretize \eqref{eq:ode} with time-step $𝚫t>0$, AD-transform the
discretized system (using all four modes under consideration, using an
approximation for the inverse of a near-identity matrix), and take the limit as
$𝚫t → 0$, yielding ODEs for the four modes.
We use $x_k = x(T₀+k𝚫t)$, so the Euler approximation
\begin{math}
  x(t+𝚫t)=x(t)+𝚫t\ddt{x}(t)
\end{math}
becomes
\begin{math}
  x_k = f_k(x_{k-1})
\end{math}
where $f_k(x) = x + 𝚫t∙g(x)$.
This makes for step-wise Jacobians
\begin{equation}
  𝐉_k = 𝐉_{f_k(x_k)} = 𝐈 + 𝚫t∙𝐉_{g(x_k)}
\end{equation}
which determine the steps of the four AD accumulation modes under
consideration,
\begin{subequations}
  \begin{align}
    \fv{x}_{k+1} &= 𝐉_k \, \fv{x}_k = \fv{x}_k + 𝚫t∙(\Jforward \; g \; x_k \; \fv{x}_k)
    \\
    \rv{x}_k &= \tran{𝐉_k} \rv{x}_{k+1} = \rv{x}_{k+1} + 𝚫t∙(\Jreverse \; g \; x_k \; \rv{x}_{k+1})
    \\
    \fiv{x}_k &= \inv{𝐉_k} \fiv{x}_{k+1}
    = \fiv{x}_{k+1} - 𝚫t∙(\Jforward \; g \; x_k \; \fiv{x}_{k+1})
    \\
    \riv{x}_{k+1} &= \itran{𝐉_k} \riv{x}_k
    = \riv{x}_k - 𝚫t∙(\Jreverse \; g \; x_k \; \riv{x}_k)
  \end{align}
\end{subequations}
using the identity
\begin{math}
  \inv{(𝐈 + 𝚫t 𝐀)} = 𝐈 - 𝚫t 𝐀 + O(𝚫t²)
\end{math}
gives per-step errors of $O(𝚫t²)$, which over the course of $O(1/𝚫t)$
time-steps gives a final numeric error of $O(𝚫t)$.
Putting these in the form $(v_{k+1}-v_k)/𝚫t = b$ and taking the limit $𝚫t→0$,
\begin{subequations}
  \begin{align}
    \ddt \fv{x} &= \Jforward \; g \; x \; \fv{x}
    & \fv{x}(T_0) &↦ \fv{x}(T_1)
    \label{eq:ddt-f}
    \\
    \ddt \rv{x} &= - \Jreverse \; g \; x \; \rv{x}
    & \rv{x}(T_0) &↤ \rv{x}(T_1)
    \label{eq:ddt-r}
    \\
    \ddt \fiv{x} &= \Jforward \; g \; x \; \fiv{x}
    & \fiv{x}(T_0) &↤ \fiv{x}(T_1)
    \label{eq:ddt-ri}
    \\
    \ddt \riv{x} &= - \Jreverse \; g \; x \; \riv{x}
    & \riv{x}(T_0) &↦ \riv{x}(T_1)
    \label{eq:ddt-fi}
  \end{align}
\end{subequations}
Note that \eqref{eq:ddt-f} and \eqref{eq:ddt-ri} are identical, as are
\eqref{eq:ddt-r} and \eqref{eq:ddt-fi}, except that the specified/calculated
boundary condition differs in location between $T₀$ and $T₁$, so the direction
of integration is reversed.
If the primal equation \eqref{eq:ode} is stable, the linear ODEs
\eqref{eq:ddt-f} and \eqref{eq:ddt-r} are also stable, and therefore
\eqref{eq:ddt-fi} and \eqref{eq:ddt-ri} would be unstable.
This is an intrinsic property: if a linear operator is stable its inverse will
be unstable, since the eigenvalues are inverted.

In a system which allows AD transforms of basis functions to be
user-specified, and which has higher-order functions including differential
equation solvers, this would suggest efficient direct transforms of such
solvers not just for forward and reverse AD, as is now routine (see for example
the Diffrax\footnote{\url{https://docs.kidger.site/diffrax/}} subsystem
\citep{Kidger-2021a} for the AD-enabled language JAX, or
torchode\footnote{\url{https://github.com/martenlienen/torchode}} for PyTorch
\citep{Lienen-Gunnemann-2022a}) but also for inverse AD.

\section{Implementation}

We have developed a prototype implementation for inverse AD of constant-width
computations.
Work is currently underway to automatically detect ``lumps'' and handle them
appropriately.
This problem is more difficult than it might initially appear, because the
computation graph can be broken up in different ways.
To avoid treating this as a brute-force combinatorial problem, either
heuristics must be employed, or connections to efficient graph algorithms like
max flow must be made.
See \S\ref{sec:constWidth} for further discussion of this issue.

\section{Related Work}

\subsection{Classic work on the inverse problem.}

The idea of direct calculation of the solution of a linear system
resulting from the linearization of a function represented as a computer program
was introduced by \citet{Griewank1990DCo} and
elaborated by \citet{Dixon1991UoA}, \citet{Utke1996ENS}, and
\cite[Chapter~4]{Hossain1998OtC}, using a framework in which the
multiple $𝐉_t$ matrices here are replaced by a single much larger
matrix.
That framework is quite general, but requires that the computation graph be
stored and manipulated in a fashion which seems difficult to migrate to
compile-time.
At root this is because that formulation is not compositional.
The present framework, which is compositional, is amenable to efficient
implementation, which we have done in a preliminary implementation.

\subsection{Recent related work.}

An early version of this work was publicly discussed in 2019\footnote{
\url{https://openreview.net/forum?id=Bygj2Ys6IS} \quad
\url{{https://openreview.net/pdf?id=Bygj2Ys6IS}}}
and its implementation in JAX was proposed in 2022 by Neil Girdhar and discussed at length,\footnote{\url{https://github.com/jax-ml/jax/issues/12494}} and added and removed from JAX proper by Matt Johnson.\footnote{\url{https://github.com/jax-ml/jax/commit/902fc0c3d2b3ec9b6034c66074984386ec35606f}}
\citet{naumann-2024a} discussed some of the ideas presented
in that work and here, but did not cite the earlier work just discussed.
It also does not contain equations \eqref{eq:inverseAD}, \eqref{eq:primitives},
or \eqref{eq:blocks}, does not discuss the fact that simple operations that
preserve width result in $𝐀:l×l$ and $𝐁:l×(k-l)$, does not discuss the fact
that although $𝐉_t$ is not structurally symmetric, $\inv{𝐉_t}$ has the same
structural sparsity as $𝐉_t$ and therefore the amount of arithmetic is the same
as for conventional Forward and Reverse modes, but that these are transposed so
Forward Inverse Mode writes to the derivative-related quantities associated
with \emph{all} involved variables of each basis function invocation, while
Reverse Inverse Mode writes only to the quantities associated with slots
\emph{written to} in the primal computation of each basis function, does not
discuss ``lumpification,'' does not show a compositional method analogous to
Figs.~\ref{fig:graphAtomicXforms}, \ref{fig:example}, and~\ref{fig:derived},
does not discuss the invariants \eqref{eq:invariants}, and does not discuss
issues of stability of inverse AD, \eg\ of inverse AD of a stable ODE.

\citet{naumann-2024a} does explore a number of technical issues not discussed
here.
It introduces a formulation in which ``pass through'' nodes are added to the
computation graph, so that width is defined in terms of node cuts rather than
edge cuts.
It also shows that a variety of scheduling issues associated with inverse AD
are NP-complete.

The NP-hardness result might not impact the formulation here as it seeks to
determine the minimal number of arithmetic operations to perform inverse AD\@.
In the case with basis functions with multiple outputs, as would arise when
coalescing lumps, this would involve selecting a lumpification that made $𝐀$
and $𝐁$ suitably sparse to minimize the number of arithmetic operations to
compute $\inv{𝐀}$ and $- \inv{𝐀} 𝐁$.
If one was not concerned with the sparsity of these operations, treating them
as atomic, and one was interested only in doing forward inverse AD or reverse
inverse AD to compute inverse Jacobian (transpose) vector products, and not
some hybrid that involved vertex, edge, or face  elimination to compute the
full inverse Jacobian, the requisite lumpification process reduces to selecting
a topological sort that minimizes maximal lump size, maximal lump width, or
maximal values of $l$ and $k$.
We conjecture that a greedy algorithm may suffice for this.

\subsection{Linguistic Support and Program Inversion.}

Invertible computation in general, and automatic program inversion in
particular, has been a subject of study in programming language theory for
decades \citep{Gries1981}, with applications in protocol design, automatic
model updates, \etc{}
This is closely related to a special case of the present formulation, namely
fully invertible numeric computations.
The idea of program inversion has been pursued by a small but dedicated
community, with some striking results \citep{Matsuda-etal-2010a,
  Matsuda-Wang-2020a, kristensen-etal-2022a}.

Here we have focused on a computation which is locally invertible, meaning that
its linearization around a point is invertible.
This property necessarily holds when a computation is globally invertible and
also differentiable.
Techniques developed there for statically guaranteeing complete invertibility
using a type system may be be applicable here as well, to guarantee local
invertibility: guaranteeing by construction that the data flow graph of active
variables will be constant width.

We also note that it is well known that if the primal is stepwise invertible
then the tape normally used in reverse mode to store values computed during the
forward sweep for use in the reverse sweep in reverse order can be eliminated
as these values can be recomputed in reverse order during the reverse sweep \citep{PEARLMUTTER94C, maclaurin2015gradient}.
Analogously, under this constraint the tape can also be eliminated from reverse-inverse mode.

\subsection{Quantum Computing and Machine Learning.}

Quantum computation is also invertible and differentiable in the precise sense
used in this manuscript.
On a topical note, invertible differentiable programs form the kernels of a
variety of generative AI systems: methods like BS-Infomax
\citep{BELL-SEJNOWSKI95A} and its context-sensitive generalizations like cICA
\citep{PARRA-ETAL95A, PEARLMUTTER-PARRA96A}, monotonic neural networks
\citep{Wehenkel-Louppe-2019a}, normalizing flows
\citep{Papamakarios-etal-2021a}, and stable diffusion
\citep{Rombach_2022_CVPR}, have at their hearts multi-layer structures
carefully constructed to be invertible, along with manually-derived inversion
procedures.
Bread-and-butter deep learning models like ResNet
\citep{behrmann2019invertibleresidualnetworks} or CNNs frequently used in deep
learning can be easily modified to maintain a constant width, particularly if
re-cast as differential equations \citep{MALEKI-ETAL-2021A}.
Efficient support for locally invertible numeric computations, using inverse
AD, would allow models of the above classes to be generalized and much more
easily implemented, as the extremely tedious and error-prone manual derivations
and coding would be avoided.

\section*{Acknowledgments}

This was supported, in part, by US National Science Foundation (NSF) grants
1522954-IIS and 1734938-IIS, by a US Intelligence Advanced Research Projects
Activity (IARPA) grant via Department of Interior/Interior Business Center
(DOI/IBC) contract number D17PC00341, by Science Foundation Ireland under grant
number 20/FFP-P/8853, and by the Defense Advance Research Projects Agency
(prime contract award HR0011222003, subcontract award 2103299-01, grant
13001129).
The content of the information does not necessarily reflect the position of the
US Government.
No official endorsement should be inferred.
Approved for public release; distribution is unlimited.

\clearpage

\bibliographystyle{siam}
\bibliography{QobiTeX,ad,abb-full,boltzmann}
\end{document}